\documentclass[11pt]{article}
\pdfoutput=1

\usepackage{a4wide}
\usepackage{Packages}

\usepackage{Definitions}

\usepackage{comment}
\usepackage{setspace}

\hypersetup{pdfauthor={D.Orlando, S.Reffert},pdftitle={Combinatorics of the Dimer Model on a Strip}}

\author{
  \begin{minipage}{.97\linewidth}
    \vspace{1cm}
    \begin{center}
      \begin{small}
\textbf{Domenico Orlando}${}^1$ and \textbf{Susanne Reffert}${}^2$
      \end{small}
    \end{center}
    \vspace{1cm}
    \hspace{2cm}\begin{minipage}{.7\linewidth}
     {\it \begin{footnotesize}
        \begin{itemize}
                \item[${}^1$] Universit\`a di Milano-Bicocca and INFN, Sezione di
          Milano-Bicocca,  \\
          P.zza della Scienza, 3, I-20126 Milano, Italy
        \item[${}^2$] Institute for Theoretical Physics, University of Amsterdam,\\
          Valckenierstraat 65, 1018XE Amsterdam, The Netherlands.
        \end{itemize}
     \end{footnotesize}}
    \end{minipage}
    \vspace{1cm}
  \end{minipage}
}

\date{}

\title{\vspace{1.5cm}
  \begin{huge}
    \textbf{Combinatorics of the Dimer Model on a Strip}
  \end{huge}
}

\begin{document}

\begin{titlepage}
  \maketitle
  \thispagestyle{empty}

  \vspace{-12cm}
  \begin{flushright}
    ITFA-2007-41
  \end{flushright}

  \vspace{14cm}

  \begin{center}
    \textsc{Abstract}\\
  \end{center}
  In this note, we give a closed formula for the partition function of
  the dimer model living on a $2\times n$ strip of squares or hexagons
  on the torus for arbitrary even $n$. The result is derived in two
  ways, by using a Potts model like description for the dimers, and
  via a recursion relation that was obtained from a map to a $1D$ monomer--dimer system.
 
  The problem of finding the number of perfect matchings can also be
  translated to the problem of finding a minmal feedback arc set on
  the dual graph.
  
 \end{titlepage}

\onehalfspace


\newpage

\section{Introduction}
\label{sec:intro}

In this note, we give a closed formula for the partition function of
the dimer model living on a $2\times n$ strip of squares or hexagons
on the torus for arbitrary even $n$.

The dimer model is concerned with the statistical mechanics of close
packed dimer arrangements on a bipartite graph.  The real--world
representation of the dimer model is the adsorption of diatomic
molecules on a crystal surface.

In the 1960s, the question of how many perfect matchings exist on a
plane graph was solved independently by Kasteleyn~\cite{Kasteleyn1,
  Kasteleyn2}, and Temperley and Fisher~\cite{Temperley, Fisher}: the
total number is given by the Pfaffian of a signed, weighted adjacency
matrix of the graph (the Kasteleyn matrix).  Much of the original
interest in the dimer model arose because it provides a simple and
elegant solution for the 2--dimensional Ising model~\cite{Hurst}.

The problem of enumerating perfect matchings is of course a classical
problem in graph theory and combinatorics (see
\emph{e.g.}~\cite{Kenyon2}), and can also be phrased in terms of
domino tilings~\cite{Kenyon1}.  During the last years, the interest in
the dimer model was revived thanks to its manifold connections to other
branches of mathematics and physics, such as the topological string
A--model~\cite{Okounkov:2003sp, Iqbal:2003ds}, real algebraic
geometry~\cite{okounkov1,okounkov2}, BPS black holes from $D$--branes
wrapping collapsed cycles~\cite{Heckman:2006sk} and supersymmetric
quantum mechanics and categorification
techniques~\cite{Dijkgraaf:2007yr}.  Furthermore, a correspondence
between the dimer model and quiver gauge theories arising from
$D3$--branes probing a singular toric surface was discovered and
worked out in great detail~\cite{Hanany:2005ve, Franco:2005rj,
  Franco:2005sm, Hanany:2005ss, Hanany:2006nm}. An explanation of this
correspondence via mirror symmetry was given in~\cite{Feng:2005gw}.

\bigskip The plan of this note is as follows. We briefly introduce the
dimer model and give a Potts--like description of the dimer model
living on a $2\times n$ strip on the torus. Using this description, we
derive a closed formula for the Newton polynomial for any value of
$n$. The same result can also be derived with a recursion relation
obtained by mapping the problem to a one--dimensional monomer--dimer 
system, and is given both for a strip of squares and a strip of
hexagons.

Furthermore, the question is translated to the problem of finding a minimal feedback arc set on the dual graph.

\bigskip

A \emph{bipartite} graph $\gra$ is a graph in which all
vertices can be colored black or white, such that each black vertex
is only connected by links to white vertices and vice versa. Let $M$
be a subset of the set $E$ of edges of $\gra$. $M$ is called a
\emph{matching}, if its elements are links and no two of them are
adjacent. If every vertex of $\gra$ is saturated under $M$, the
matching is called \emph{perfect}. Such a link that joins a black and
a white vertex is called a \emph{dimer}. The \emph{dimer model} describes the
statistical mechanics of a system of random perfect matchings.  In the
simplest case, we ask for the number of close packed dimer
configurations, \emph{i.e.} the number of perfect matchings.

Kasteleyn~\cite{Kasteleyn1, Kasteleyn2} introduced an orientation on
\gra, which leads to a signed adjacency matrix $K$, now called the
\emph{Kasteleyn matrix}. The Pfaffian of $K$ gives the number of
perfect matchings.
A \emph{Kasteleyn orientation} fulfills the following
condition: the product of all edge weights around a face must equal $-1$
if the number of edges around the face is $0 \mod 4$.  If the number
of edges equals $2 \mod 4$, the product must equal $1$~\cite{Kenyon2}. One
can choose an orientation by consistently assigning arrows to
the edges of the graph, as originally suggested by Kasteleyn
\cite{Kasteleyn1, Kasteleyn2}. 
The above treatment can be straightforwardly generalized to any genus
$g$ Riemann surface.  

\begin{figure}
  \begin{center}
    \includegraphics[width=45mm]{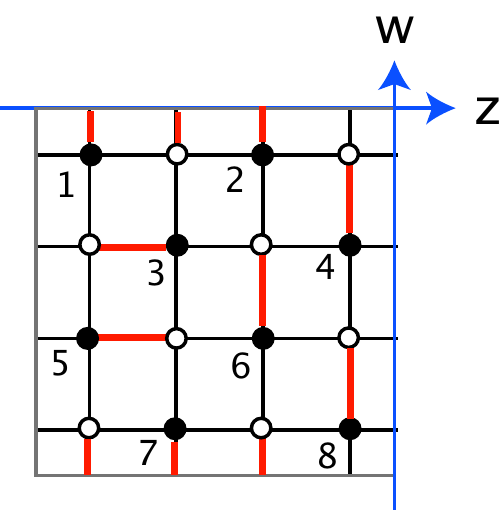}
    \caption{Example of a square graph on the torus}
    \label{fig:example}
  \end{center}
\end{figure}

In the following, we will restrict ourselves to regular $2\times n$ graphs ${\cal G}_{2,n}$ embedded on a 
torus, to which we will refer in the following as a \emph{strip}. 
On the torus, there are two
non--trivial cycles, which we will denote by $z$ and $w$. \\
In the case
of the plane graph, the edge weights originated solely from the
Kasteleyn orientation. We choose a positive direction on the dimers,
say $\bullet\to\circ$. Now we assign the weight $z$ ($w$) to each edge
which crosses the cycle $z$ ($w$) in positive direction and the weight
$1/z$ ($1/w$) to each edge which crosses it in negative
direction. While the Pfaffian of the Kasteleyn matrix yielded a number
in the case of the plane graph, it becomes a polynomial in $z$ and
$w$ on the torus, the so--called characteristic polynomial or \emph{Newton
polynomial} of the graph. The coefficient of each monomial $z^pw^q$
gives the number of matchings with \emph{weight} $(z,w)=(p,q)$.
These are matchings with the number of dimers crossing $z$ in positive
direction minus the number of dimers crossing $z$ in negative
direction equal to $p$ (analogous for $q$). In the literature,
what we call the weight is usually referred to as the slope of a height
function defined on the composition of two matchings.\footnote{The height
function is defined as follows. Choose a reference matching $\mathrm{PM}_0$. To
find the slope of a matching $\mathrm{PM}$, compose it with the reference
matching, $\mathrm{PM}-\mathrm{PM}_0$, where the minus serves to change the orientation of
$\mathrm{PM}_0$ to $\circ\to\bullet$. This results in closed loops (composition
cycles) and double line dimers. The rule is that when an edge in $\mathrm{PM}$
belonging to a closed loop is crossed such that the black node is to
its left (right), the height changes by $+1$ $(-1)$. If an edge
belonging to $\mathrm{PM}_0$ is crossed, the signs are reversed. This height
function is defined up to the choice of the reference matching
$\mathrm{PM}_0$. Crossing the boundary of the fundamental region of the torus, this function
can jump. If the height function jumps by $p$ units crossing $z$, it
is associated to the power $z^p$ in the Newton polynomial of the graph
(and equivalently for $w$). Choosing a different reference matching
results in a common prefactor of $z^{p_0}w^{q_0}$ for all
monomials. Our method of assigning weights to a matching corresponds
to choosing a reference matching of weight $(0,0)$ that does not
intersect the $z$ or $w$ cycle.} 
The matching shown in Figure
\ref{fig:example} has weight $(1,0)$, where $1=2-1$.

The partition function, or Newton polynomial, of the dimer model on
the torus takes the form
\begin{equation}
\label{eq:torus-partition}
  {\mathcal{P}}_{m,n}(z,w)=\sqrt{\det\,K}=\sum_{n_z, n_w} N_{n_z,n_w}\,(-1)^{n_z+n_w+n_zn_w}z^{n_z}w^{n_w},
\end{equation}
where the $N_{n_z,n_w}$ count the number of matchings of weight
(height change) $(n_z,n_w)$.
Furthermore, the total number of matchings for the square graph on the torus is given by
\begin{equation}\label{eq:totalnumber}
Z_{m,n}=\frac{1}{2}\left(-{\mathcal{P}}_{m,n}^{\,\mathrm{sq}}(1,1)+{\mathcal{P}}_{m,n}^{\,\mathrm{sq}}(1,-1)+{\mathcal{P}}_{m,n}^{\,\mathrm{sq}}(-1,1)+{\mathcal{P}}_{m,n}^{\,\mathrm{sq}}(-1,-1)\right),
\end{equation}
where the first term is always zero.


\section{Combinatorics for the strip of squares}
\label{sec:strip}

\begin{figure}
  \begin{center}
    \includegraphics[width=70mm]{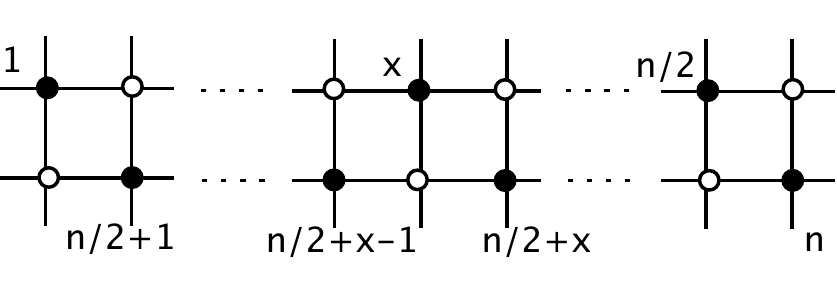}
    \caption{Long strip of squares on the torus}
    \label{fig:strip}
  \end{center}
\end{figure}

We consider a long strip of squares $\gra_{2,n}$ on the torus, see
Figure~\ref{fig:strip}, focusing on the $z$--weights.
A strip containing $n$ black nodes can accommodate matchings with
weights $n/2, n/2-1, \dots, -(n/2-1),-n/2$, \emph{i.e.} there are $n+1$
subsets. We would like to find a direct way of obtaining the
multiplicities of the matchings of a given weight, \emph{i.e.} the
information contained in the characteristic polynomial, with just the
number $n$ as input data. Using the Kasteleyn construction, we find
the following multiplicities for the first five cases:
\begin{equation}
  \label{eq:seq}
  \begin{array}{ccccccccccccc}
    n && z^{-5} & z^{-4} & z^{-3} & z^{-2} & z^{-1} & z^{0} & z^{1} & z^{2} & z^{3} & z^{4} & z^{5} \\ \hline
    2 &&&&&& 1 &4 & 1 \\
    4 &&&&& 1 &8 &16&8& 1 \\
    6 &&&& 1 & 12 & 48 & 76 &48 & 12& 1 \\
    8 &&& 1&16&96&272&384&272&96&16&1 \\
    10&& 1&20&160&660&1520&2004&1520&660&160&20&1
  \end{array}
\end{equation}
The above sequences do not have an obvious structure.
We will solve the problem using an operator perspective.

\subsection{Potts picture derivation}

In order to find a convenient formalism, we use the following picture:
we attach to each black node a $\setZ_m$ spin, which can point along
all the directions in which the black node is joined to a white node
by an edge. This results in a description reminiscent of the
$m$--state Potts model, where $m$ here is the valency of the nodes. A
hexagon graph results in a 3--state model, while the square graph
gives a 4--state model, in which the spin can point up, down, left, or
right: \up, \down, \rright, \lleft. To describe a dimer configuration,
we take the spins to point to those white nodes which are joined to
the black nodes by a dimer. Since we are interested in perfect
matchings, we must restrict the possible configurations of the Potts
model to those, in which each node is only touched by one dimer.

We label the black nodes by numbers as shown in Figure
\ref{fig:example}. The perfect matching shown in Figure
\ref{fig:example} can be written as the following spin state:
\begin{equation}
  \label{eq:examplespin}
  \ket{ \ua\ua\leftarrow\ua\ra\ua\da\ua} \, .
\end{equation}

As we will show in the following, it is possible to define basic operations on the spins allowing us to reach all states starting from the highest weight state.

\begin{figure}
  \begin{center}
    \begin{minipage}{.85\linewidth}
      \subfigure[$\mathrm{d}_a(x)$]{\includegraphics[width=.23\textwidth]{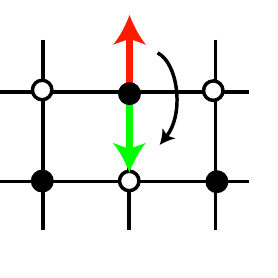}}
      \hfill
      \subfigure[$\mathrm{d}_b(x, n/2+x-1)$]{\includegraphics[width=.23\textwidth]{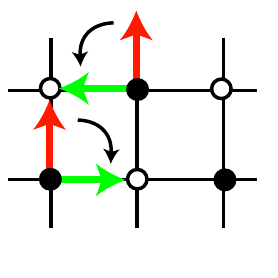}}
      \hfill
      \subfigure[$\mathrm{d}_b(x, n/2+x)$]{\includegraphics[width=.23\textwidth]{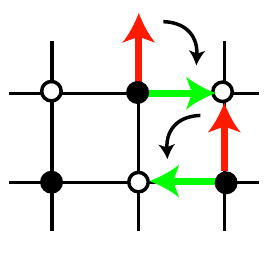}}
    \end{minipage}
    \caption{The three basic operations that can be performed on a spin pointing up}
    \label{fig:operations}
  \end{center}
\end{figure}

The highest weight state is unique and in the spin picture is the one with all spins up:
\begin{equation}\label{eq:highest}
 \ket{\underbrace{\ua\ua\ua\dots\ua\ua}_{n} } =: \ket{n} \, .
\end{equation}
Also the lowest weight state is unique and is the one with all spins down. We are now looking for basic operations on the spins, which take us from the highest weight state to another state. A spin pointing up can be flipped to three possibilities: down, left and right. Given the structure of the highest weight matching, flipping a spin down gives another perfect matching, while flipping the spin on node $x$ left requires us to flip the spin on node $n/2+x-1$ right to arrive at a perfect matching. Flipping the spin on node $x$ right requires us to flip the spin at $n/2+x$ to the left. See Figure~\ref{fig:operations} for these three basic operations.
Each of these operations lowers the weight of the perfect matching by
one.  We shall denote the spin flip from up to down at node $x$ by
${\di}_a(x)$.  The spin flips by $\pi /2$ on the nodes $x$ and $y$
will be denoted by $\di_b(x,y)$.  For $x = 2 , \dots , n/2 $, $y=
n/2+x-1$ or $y=n/2+x$ and for $x=1$, $y=n/2+x$ or $y=n$.

To familiarize ourselves with this operator picture, see the example of $n=4$ black nodes in Appendix~\ref{sec:ex}.

\begin{figure}
  \begin{center}
    \includegraphics[width=.55\textwidth]{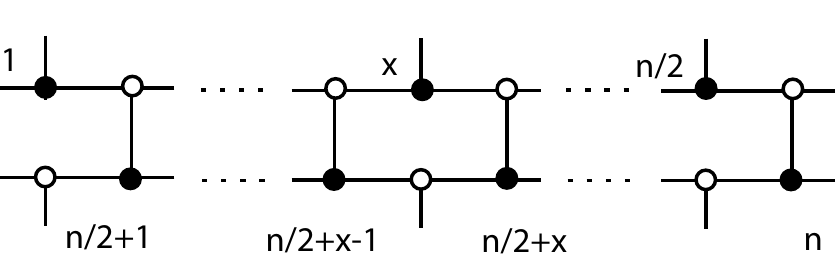}
    \caption{Long strip of hexagons on the torus}
    \label{fig:hex_strip}
  \end{center}
\end{figure}

Before we treat the square graph, we first solve the case of the strip of hexagons, which can be obtained from the square by deleting half of the vertical links, see Figure~\ref{fig:hex_strip}.
\subsubsection{Combinatorics for the strip of hexagons}
\label{sec:hex_strip}
%
Since the strip of hexagons is a 3--state model and does not accommodate downwards (or upwards,
depending on which half of the vertical links we choose to delete)
pointing spins, the weights in $z$ only
run from 0 to $n/2$ and only the spin flips by $\pi /2$, $\di_b(x,y)$,
exist.

The first five cases have the following multiplicities:
\begin{equation}
  \begin{array}{ccccccc}
    n & z^0 & z^1 & z^2 & z^3 & z^4 & z^5\\ \hline
    2 &2 & 1 \\
    4 &2 &4 &1 \\
    6 &2 &9 & 6& 1 \\
    8 &2 &16&20& 8 &1 \\
    10&2 &25&50&35 &10 & 1
  \end{array}
\end{equation}
The combinatorics therefore stems exclusively from repeated operations
of $\di_b$ on the state $\ket{n}$.  On the hexagon graph, the example of $n=4$ black nodes
shown in the Appendix reduces to what is represented in
Figure~\ref{fig:hex_sequence}.
\begin{figure}
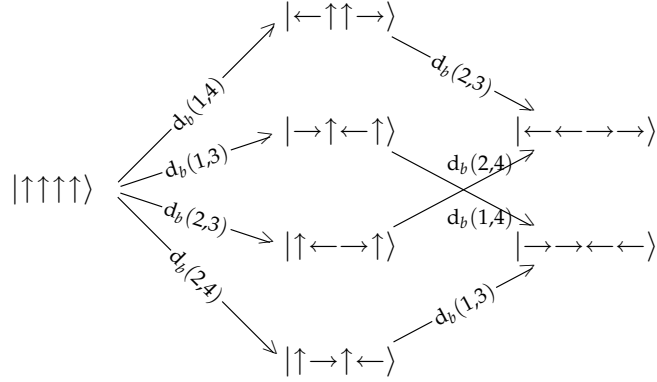

  \begin{center}
\begin{diagram}
  &&&&&&&&&&&&\\
  &&&&& \ket{\la \ua \ua \ra} \quad&&&&&&&&\\
   &&&\ruTo(3,3)~{\di_b(1,4)}&&&\rdTo(3,2)~{\di_b(2,3)}&&&&&&\\
  &&&&& \ket{\ra \ua \la \ua} \quad&&&\quad \ket{\la\la\ra\ra}&&&&&\\
 \ket{\ua \ua \ua \ua}&&&  \ruTo(3,1)~{\di_b(1,3)}&&&\rdTo(3,2)^{\di_b(2,4)}& \ruTo(3,2)_{\di_b(1,4)}&&&&&&&& \\
  &&\rdTo(3,1)~{\di_b(2,3)} \rdTo(3,3)~{\di_b(2,4)} &&& \ket{\ua \la \ra \ua} \quad&&&\quad \ket{\ra \ra \la \la}&&&&&\\
  &&&&&&&\ruTo(3,2)~{\di_b(1,3)}&&&& \\
  &&&&& \ket{\ua \ra \ua \la}\quad &&&&&&\\
\end{diagram}
     \caption{Strip of hexagons with $n=4$ black nodes}
    \label{fig:hex_sequence}
  \end{center}
\end{figure}
\begin{theo}
The general formula for the number of perfect matchings with $z$--weight $n/2-p$ on the strip of hexagons on the torus is
\begin{equation}
  \label{eq:anp_hex}
  a_{n,p}=\frac{n}{p!}\prod_{q=1}^{p-1}\,(n-p-q) = \frac{n}{n - p} \binom{ n - p}{p}\,.
\end{equation}
\end{theo}

\begin{proof}
  The operator description explained above is convenient because it allows to map the
  two--dimensional dimer problem on the hexagonal strip to a
  one--dimensional monomer--dimer problem. In fact we can simply
  concentrate on the upper line of the strip and consider two nodes to
  be occupied by a dimer if they are occupied by a horizontal line in the strip resulting from a $b$--move, and to be occupied by a
  monomer if the dimer is vertical in the strip. The occupation of the lower line is uniquely determined
  by consistency. Let us now write down the
  partition function for a monomer--dimer system living on a $1D$ lattice of
  $n$ nodes with periodic boundary conditions (note that we do not require the lattice to be bipartite anymore, so $n$ can be odd). The
  partition function reads:
  \begin{equation}
    \label{eq:12}
    Q_n (q) = \sum_{k=0}^{\lfloor n/2 \rfloor} a_{n,p}\, q^p \, ,
  \end{equation}
  where $a_{n,p}$ is the number of configurations with $p$ dimers.
  Let us first start with a slightly simpler system of $n$ nodes on a
  line with free boundary conditions, and let $P_n (q)$ be the
  corresponding monomer--dimer partition function. This function satisfies the following recursion relation:
  \begin{equation}
    P_{n+1} (q) = P_{n} (q) + q\, P_{n-1} (q) \, .  
  \end{equation}
  This can be understood as follows. When adding an extra point after
  the $n$-th, one can either add a monomer which leaves $P_{n+1} (q) =
  P_n (q)$ (since $P_n(q)$ counts the dimers), or add a dimer (and
  multiply by $q$) if the last point was previously occupied by a
  monomer. The configurations in $P_n(q)$ where the $n$-th point is a
  monomer are precisely counted by $P_{n-1} (q)$. Note that this is
  the $q$--analogue of the Fibonacci sequence and in fact the actual
  Fibonacci sequence can be recovered by $\set{P_n(1)}$, where the
  $P_n(q)$ are obtained by using the initial conditions $P_1(q) = 1$,
  $P_2(q) = 1 + q$.

  In a similar way, we can understand the relation between $P_n(q) $
  and $Q_n(q)$. In fact, making the plane strip periodic by adding a
  line between the first and the $n$-th node can either leave the
  partition function invariant or add an extra dimer for each
  configuration where the first and last nodes are occupied by a
  monomer. Therefore,
  \begin{equation}
    Q_n (q) = P_n (q) + q\, P_{n-2} (q) \, .    
  \end{equation}
  It is easy to see that $Q_n (q)$ satisfies the same recursion
  relation as $P_n(q)$. On the other hand, we are interested in the
  $n$ even case, so it is better to recast it in the form
  \begin{equation}\label{eq:recursionQ}
    Q_{n+2} (q) = \left( 1 + 2q \right) Q_{n} (q) - q^2  Q_{n-2} (q) \, ,
  \end{equation}
  that can be easily solved with the initial conditions
  \begin{align}
    Q_0(q) = 1 \, , && Q_2(q) = 1 + 2\, q \, ,
  \end{align}
  the solution being
  \begin{equation}\label{eq:Q_explicit}
    Q_n (q) = \frac{1}{2^n} \left[ \left( 1 - \sqrt{1 + 4\, q} \right)^n  + \left( 1 + \sqrt{ 1 + 4\, q} \right)^n\right] \, . 
  \end{equation}
 To expand~(\ref{eq:Q_explicit}) in powers of $q$ we can use the identity
  \begin{equation}
    \left( a + b \right)^m+ \left( a - b \right)^m = 2 \sum_{k=0}^{\lfloor m/2 \rfloor}{m \choose 2\,k}\,a^{m-2k}\,b^{2k} \, ,
  \end{equation}
  and find
  \begin{equation}
    Q_n (q) = \frac{2}{2^n} \sum_{k=0}^{\lfloor n/2 \rfloor} \binom{n}{2\,k} \left( 1 + 4\, q \right)^k =  \frac{1}{2^{n-1}} \sum_{k=0}^{\lfloor n/2 \rfloor} \binom{n}{2\,k} \sum_{p=0}^k \binom{k}{p} \left(4\, q \right)^p \, . 
  \end{equation}
  Inverting the sums, we find
  \begin{equation}
    Q_n (q) = \frac{1}{2^{n-1}} \sum_{p=0}^{\lfloor n/2 \rfloor} \left(4\, q \right)^p \sum_{k=p}^{\lfloor n/2 \rfloor}   \binom{n}{2\,k} \binom{k}{p} \,.
  \end{equation}
  Using the identity (see~\cite{Riordan:1968})
  \begin{equation}
    \sum_{k=p}^{\lfloor n/2 \rfloor}   \binom{n}{2k} \binom{k}{p}  = 2^{n-1-2p} \frac{n}{n-p} \binom{n-p}{p},
  \end{equation}
  we obtain
  \begin{equation}
    Q_n (q) = \sum_{p=0}^{\lfloor n / 2 \rfloor} \frac{n}{n-p} \binom{n-p}{p} q^p \, .\qedhere   
  \end{equation}
\end{proof}
  The explicit expression for $Q_n(q)$ allows to compactly summarize
  the result for any strip by defining the generating function
  \begin{equation}
    \tilde {\mathcal{F}}^{\text{hex}} (s, q) = \sum_{n=0}^\infty s^n Q_n (q) = \frac{ s - 2 }{s^2 q + s - 1 } \, .
  \end{equation}

\begin{corr}
  The partition function for the strip of hexagons on the torus,
  considering only $z$ is:
  \begin{multline}
    \label{eq:poly_hex}
    {\mathcal{P}}_n^{\,\mathrm{hex}}(z) = z^{n/2}\,\sum_{p=0}^n(-1)^p\,z^{-p}\,a_{n,p} = z^{n/2} Q_n (-1/z) = \\
    = \frac{z^{n/2}}{2^n} \left[ \left( 1 - \sqrt{1 - \frac{4}{z}} \right)^n  + \left( 1 + \sqrt{ 1 - \frac{4}{z}} \right)^n\right] \, . 
  \end{multline}
\end{corr}
Since there is only one matching with weight $w$ and one with weight
$1/w$, the full Newton polynomial on the torus is
\begin{equation}
  \label{eq:poly_full_hex}
  {\mathcal{P}}_{2,n}^{\,\mathrm{hex}}(z,w)={\mathcal{P}}_n^{\,\mathrm{hex}}(z)-w-\frac{1}{w}.
\end{equation}

Again we can summarize the result into a generating function in $z$ for the dimers on a hexagonal strip as follows:
\begin{equation}
  \mathcal{F}^{\text{hex}} (s, z) = \sum_{n=0}^\infty {\mathcal{P}}_n^{\,\mathrm{hex}}(z) s^n = \frac{ 2 - \imath s \sqrt{z} }{1 - \imath p \sqrt{z} - s^2 } \, . 
\end{equation}

\subsubsection{Generalization to the square strip}

\begin{figure}[h!]
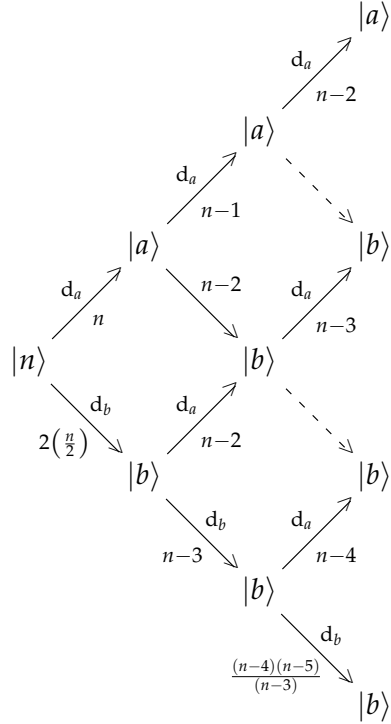

  \begin{center}
\begin{diagram}
  &&&&&& \ket{a} \\
  &&&&& \ruTo^{\di_a}_{n-2}\\
  &&&& \ket{a} \\
  &&& \ruTo^{\di_a}_{n-1} && \rdDashto\\
  && \ket{a} &&&& \ket{b} \\
  & \ruTo^{\di_a}_n && \rdTo{n-2} && \ruTo^{\di_a}_{n-3}\\
  \ket{n} & & & & \ket{b} \\
  & \rdTo^{\di_b}_{2 \left(\frac{n}{2} \right)} && \ruTo^{\di_a}_{n-2} && \rdDashto\\
  && \ket{b} &&&& \ket{b} \\
  &&& \rdTo^{\di_b}_{n-3} && \ruTo^{\di_a}_{n-4}\\
  &&&& \ket{b} \\
  &&&&& \rdTo^{\di_b}_{\frac{\left(n-4 \right)\left(n-5 \right)}{\left(n-3 \right)}}\\
  &&&&&& \ket{b}   
\end{diagram}
     \caption{Diagram of $\di_a, \, \di_b$ operations}
    \label{fig:diagram}
  \end{center}
\end{figure}

After having solved the problem for the strip of hexagons, we return
to the square strip.  We denote a state which only contains up and
down spins by $\ket{a}$, while a state which contains at least one
left/right pair by $\ket{b}$.  Under the action of our two
operators, we have
\begin{equation}
  \di_a \ket{a} = \ket{a} \, , \quad \di_a \ket{b} = \ket{b} \, , \quad \di_b \ket{a} = \ket{b} \, , \quad \di_b \ket{b} = \ket{b} \, .
\end{equation}
We denote by $p$ the number of times we act on the highest weight
state. The resulting states have weight $n/2-p$. The order in which
the state was acted on by the $\di_a$ and $\di_b$ does not matter, so
at level $p$, there are $p+1$ combinations of $\di_a$ and $\di_b$. On
the highest weight state $\ket{n}$, we can act in $n$ ways with
$\di_a$ and in $n$ ways with $\di_b$. On $\di_a \ket{n}$, we can then
act in $n-1$ ways with $\di_a$ and in $n-2$ ways with $\di_b$. On
$\di_b \ket{n}$, we can act in $n-2$ ways with $\di_a$ and in $n-3$
ways by $\di_b$, etc. It is in general easier to compute the number of
possibilities of acting with $\di_a$ on a given state than with
$\di_b$, so we choose the most convenient path to obtain the full
combinatorics, see Figure~\ref{fig:diagram}.

We summarize the results for the first three levels:
\begin{equation}
    \label{eq:table_comb}
  \begin{tabular}{cccccccc}
    $p$\\
    \multirow{2}{*}{$1$} &  $d_a$ & \multirow{2}{*}{$+$} & $d_b$ \\
     & $n$ &  & $n$ \\
    \multirow{2}{*}{$2$} & $d_a d_a$ & \multirow{2}{*}{$+$} & $d_a d_b$ & \multirow{2}{*}{$+$} & $d_b d_b$ \\
    &  $\frac{n \left(n - 1 \right)}{2}$ &  & $n \left(n-2 \right)$ &  & $\frac{n \left( n - 3 \right)}{2}$ \\
    \multirow{2}{*}{$3$} & $d_a d_a d_a$ & \multirow{2}{*}{$+$} & $d_a d_a d_b$ & \multirow{2}{*}{$+$} & $d_b d_b d_a$ & \multirow{2}{*}{$+$} & $d_b d_b d_b$ \\
    & $\frac{n \left(n-1 \right) \left( n -2 \right)}{3!}$ & & $\frac{n \left( n - 2 \right) \left( n - 3 \right)}{2!}$ & & $\frac{n \left( n - 3 \right) \left( n - 4 \right)}{2!}$ & & $\frac{n \left( n - 4 \right) \left( n - 5 \right)}{3!} $
  \end{tabular}
\end{equation}
Note the denominators, which avoid an overcounting of states obtained
by repeatedly applying the same operators. The result for the $p$-th
level is a sum of $p+1$ terms. The denominator of the $q$-th term on
level $p$ is $q! \left( p - q \right)!,\quad q=0,..,p$.
\begin{theo}
  \label{theo_sqaure}
  The $q$-th term on level $p$ equals to
  \begin{equation}
    \label{eq:anpq}
    b_{n,p,q}=\frac{1}{q!(p-q)!}\,n\,\frac{(n-q-1)!}{(n-q-p)!}=\frac{n}{n-q}\binom{p}{q}\binom{n-q}{p},\quad q = 0, \dots, p \, .
  \end{equation}
  The total number of perfect matchings with weight $n/2-p$ is
  therefore
  \begin{equation}
    \label{eq:anp}
    b_{n,p}=\sum_{q=0}^p\,\frac{n(n-q-1)!}{{q!(p-q)!}(n-q-p)!} \, .
  \end{equation}
\end{theo}
Note that $b_{n, n-p} = b_{n,p}$, which mirrors the symmetry of the
sequence.
\begin{proof}
  We will prove~(\ref{eq:anpq}) by induction. The first three terms are
  shown explicitly in~(\ref{eq:table_comb}). Note that $\left( p - q
  \right)$ counts the number of times $\di_a$ has been applied and $q$
  counts the number of times $\di_b$ has been applied. One can act
  with $\di_a$ on a given configuration at level $(p,q)$ in $n-q-p$
  ways, since of the $n$ upwards pointing arrows in the highest weight
  state, $(p-q)$ were turned down by acting with $\di_a$, whereas
  $2\,q$ were turned horizontally by acting $q$ times with
  $\di_b$. This number has to be normalized by the degeneracy factor
  $q!(p-q)!$. Indeed,
  \begin{equation}
    q!(p+1-q)! \, b_{n,p+1,q} = \left( n -  p - q \right)  q! \left( p - q \right)!\,b_{n,p,q} \, .
  \end{equation}
  Since every possible state can be arrived at by acting with $\di_a$
  on a state created by acting with only $\di_b$ on $|n\rangle$, it is
  enough to now prove~(\ref{eq:anpq}) for the case $p=q$, which
  corresponds to the case of the strip of hexagons. In fact
  $b_{n,p,p} = a_{n,p}$ which was calculated in the last section.
\end{proof}
\begin{corr}
  The Newton polynomial in $z$ for the strip with $n$ black nodes ($n$
  even) is
  \begin{equation}
    \label{eq:poly}
    {\mathcal{P}}_n^{\,\mathrm{sq}}(z) = z^{-n/2}\,\sum_{p=0}^n(-1)^p\,z^p\,b_{n,p} = z^{-n/2}\,\sum_{p=0}^n(-1)^p\,z^p\,\sum_{q=0}^p\,\frac{n(n-q-1)!}{{q!(p-q)!}(n-q-p)!}.
  \end{equation}
\end{corr}
This constitutes a compact sum formula for the Newton polynomial for
the strip of squares of arbitrary length.  Since there is only one
matching with weight $w$ and one with weight $1/w$, the full Newton
polynomial on the torus is
\begin{equation}
  \label{eq:poly_full}
  {\mathcal{P}}_{2,n}^{\,\mathrm{sq}}(z,w)={\mathcal{P}}_n^{\,\mathrm{sq}}(z)-w-\frac{1}{w}.
\end{equation}

\subsubsection{Recursion relation}
\label{sec:rec_strip}

Even though we already have the result for the partition function in the form of a sum, we will now derive a recursion relation which expresses ${\mathcal{P}}_n^{\,\mathrm{sq}}(z)$ through ${\mathcal{P}}_{n-2}^{\,\mathrm{sq}}(z)$, \emph{i.e.} we take a strip ${\cal G}_{2,n-2}$ and add in another square, resulting in ${\cal G}_{2,n}$. Like in the case of the hexagon strip, this will provide us with a closed form, and like for the hexagon, the derivation invokes the map to a monomer--dimer system.

\begin{figure}
  \centering
  \includegraphics[width=.35\textwidth]{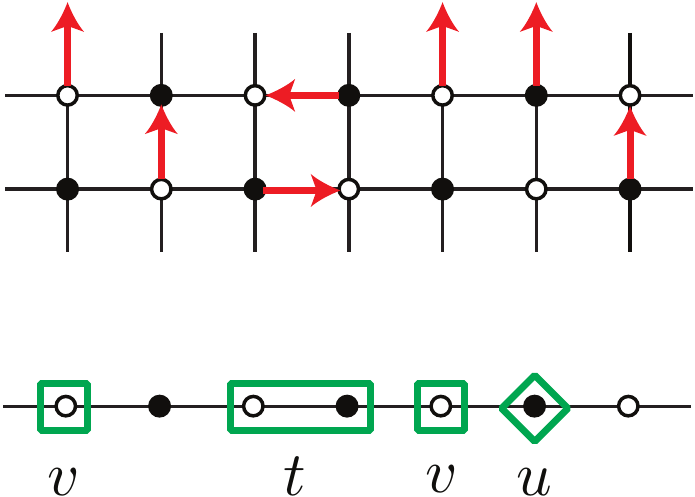}
  \caption{Mapping of a square strip configuration to a monomer--dimer
    configuration. This particular configuration contributes as $u v^2
    t$ to the partition function.}
  \label{fig:square-monomer}
\end{figure}

\begin{theo}
  The recursion relation for the partition function of the strip of
  squares on the torus is
  \begin{equation}
    \label{eq:rec_poly}
    {\mathcal{P}}_n^{\,\mathrm{sq}}(z)={\mathcal{P}}_2^{\,\mathrm{sq}}(z)\,{\mathcal{P}}_{n-2}^{\,\mathrm{sq}}(z)-{\mathcal{P}}_{n-4}^{\,\mathrm{sq}}(z)
  \end{equation}
  with the initial conditions
  \begin{align}\label{eq:init}
    \mathcal{P}_0^{\,\mathrm{sq}} (z) = 1 \, , && \mathcal{P}_2^{\,\mathrm{sq}} (z) = 4 - z - \frac{1}{z} \, .
  \end{align}
  The solution to~(\ref{eq:rec_poly}) is given by
  \begin{equation}\label{eq:sol}
    {\mathcal{P}}_n^{\,\mathrm{sq}}(z) = \frac{1}{\left( -4\, z \right)^{n/2}} \left[ \left( z - 1 - \sqrt{1 + z \left( z - 6 \right)  } \right)^n + \left( z - 1 + \sqrt{1 + z \left( z - 6 \right) } \right)^n \right]
  \end{equation}
  This result is equivalent to~(\ref{eq:poly}) and provides a closed form.
\end{theo}
\begin{proof}
  The dimer model on the square strip can be mapped to a
  one--dimensional monomer--dimer in the following way (see
  Fig.~\ref{fig:square-monomer}). A dimer living on a $z$ link
  corresponds to a $u$--monomer living on the $\bullet$ node, a dimer
  living on a $1/z$ link corresponds to a $v$--monomer on the $\circ $
  node, a dimer living on an internal vertical link corresponds to an
  empty point, and a horizontal dimer is again a dimer. In this way, the
  square strip problem becomes a monomer$^3$--dimer on a one--dimensional
  $n$--node lattice. As for the case of the hexagon graph, let us start with free boundary conditions described by the partition function
  \begin{equation}
    P_n (u,v,t) = \sum_{i,j,k} c_{i,j,k} u^i v^j t^k \, ,  
  \end{equation}
  where $c_{i,j,k} $ is the number of configurations with $i$
  $u$--monomers, $j$ $v$--monomers, $k$ dimers and $\left( n - i - j -
    2k \right)$ free nodes.

  We can always suppose without loss of generality that the first node
  is a $\bullet$, such that $P_1 (u,v,t) = 1 + u$. Adding a
  $\,\bullet\,$ node to a string of $2\,k$ points results in
  \begin{equation}
    P_{2\,k + 1 } ( u, v, t ) = \left( 1 + u \right) P_{2\,k} ( u, v, t ) + t \,P_{2\,k -1} ( u, v, t ) \, ,
  \end{equation}
  which can be read as stating that the new point is either free,
  occupied by a $u$--monomer or, if the $2\,k$-th node was free, by a new
  dimer. Similarly, adding a $\,\circ\,$ node to a string
  of $2\,k + 1$ gives
  \begin{equation}
    P_{2\,k + 2}  ( u, v, t ) = \left( 1 + v \right) P_{2\,k+1} ( u, v, t ) + t \,P_{2\,k } ( u, v, t ) \, .  
  \end{equation}
  Adding periodic boundary conditions means adding a new line between
  nodes $1$ and $n$. One can see that the corresponding partition
  function is given by
  \begin{equation}
    Q_{n} (u, v, t) = P_n (u, v, t ) + t \,P_{n-2} (u,v,t)\,,  
  \end{equation}
  which can be understood as saying that we get either the same
  configurations or a new dimer if the extremal nodes were both
  empty. Since the relation between $Q$ and $P$ is linear, they
  satisfy the same recurrence equation. On the other hand, $Q_n$ is
  well defined only for $n$ even, so the equation is better cast into the
  form
  \begin{equation}
    Q_{n+2} ( u,v, t) = \left[ \left( 1 + u \right) \left( 1 + v \right) + 2t \right] Q_n (u,v,t) - t^2 Q_{n-2} (u,v,t) \, ,  
  \end{equation}
  and solved with the initial conditions
  \begin{align}
    Q_0 (u,v,t) = 1 \, , && Q_2 (u,v,t) = 1 + u + v + u v + t \, .
  \end{align}
  
  The partition function for the dimer on the square strip is then
  obtained using the weights as they were defined on the initial
  bipartite graph and reads
  \begin{equation}
    \mathcal{P}_n^{\,\mathrm{sq}} ( z ) = Q_n ( -z, - 1/z, 1 ) \, .  
  \end{equation}
  Hence it satisfies
  \begin{equation}
    \mathcal{P}_{n+2}^{\,\mathrm{sq}} (z) = \left( 4 - z - \frac{1}{z} \right) \mathcal{P}_n^{\,\mathrm{sq}} (z) - \mathcal{P}_{n-2}^{\,\mathrm{sq}} ( z ) \, ,  
  \end{equation}
  with the initial conditions~(\ref{eq:init}).  Solved explicitly, the
  closed form~(\ref{eq:sol}) is obtained.
\end{proof}

The expression in~(\ref{eq:sol}) can be summarized by introducing a generating function in $z$ as follows.
\begin{equation}
  \mathcal{F}^{\text{sq}} (s,z) = \sum_{n=0}^\infty \mathcal{P}^{\text{sq}}_n (z) s^n = \frac{4 \imath s \left( z - 1 \right)\sqrt{z} + 8 z}{4 \imath s \left( z -1  \right)\sqrt{z} + 4 z - 4 s^2 z} \, .
\end{equation}

\subsection{Comparison with Kasteleyn's product formula}

For the total number of matchings on an $m\times n$ square lattice on
the torus, a product formula exists~\cite{Kasteleyn2}:
\begin{equation}
\label{eq:product_gen}
  Z_{m,n} = - \frac{1}{2} \sum_{a,b=0}^1 (-1)^{a+b+ab} \prod_{k=1}^{m/2}\,\prod_{l=1}^{n}2\,\sqrt{x^2\sin^2\frac{(2\,k - a  )\,\pi}{m}+x'^2\sin^2\frac{(2\,l- b)\,\pi}{n}} \, .
\end{equation}
Here, $x,\,x'$ are the edge weights which in our case are equal to one. When specializing formula~(\ref{eq:product_gen}) to $m=2$, this gives
\begin{eqnarray}\label{eq:product_strip}
  Z_{2,n}&=&-\frac{1}{2}\,\prod_{l=1}^{n}2\,\left|\,\sin\frac{2\pi l}{n}\,\right|+\frac{1}{2}\,\prod_{l=1}^{n}2\,\left|\,\sin\frac{(2l-1)\pi}{n}\,\right| + \nonumber\\
  &&+\frac{1}{2}\,\prod_{l=1}^{n}2\,\sqrt{1+\sin^2\frac{2\pi l}{n}}+\frac{1}{2}\,\prod_{l=1}^{n}2\,\sqrt{1+\sin^2\frac{(2l-1)\pi}{n}} = \nonumber\\
  &=:&-A_1+A_2+A_3+A_4.
\end{eqnarray}
Note, that the first term, $A_1$, is identical to zero (already for general $m$). Furthermore, $A_2=2$, and
\begin{equation}\label{eq:diffbytwo}
A_2-A_1=2,\quad A_4-A_3=2\quad \forall\,n.
\end{equation}
From the product formula~(\ref{eq:product_strip}), this is little
obvious, but clear from the point of view of formula
(\ref{eq:totalnumber}) and the Newton polynomial
${\mathcal{P}}_{2,n}^{\,\mathrm{sq}}$. Since the monomials in $w$ are
$-w$ and $-1/w$, ${\mathcal{P}}_{2,n}^{\,\mathrm{sq}}(\,\cdot\,, 1)$
and ${\mathcal{P}}_{2,n}^{\,\mathrm{sq}}(\,\cdot\,, -1)$ must differ
by four. Therefore, we arrive at the identifications
\begin{equation}\label{eq:identify}
A_1=\tfrac{1}{2}\,{\mathcal{P}}_{2,n}^{\,\mathrm{sq}}(1,1),\quad A_2=\tfrac{1}{2}\,{\mathcal{P}}_{2,n}^{\,\mathrm{sq}}(1,-1),\quad A_3=\tfrac{1}{2}\,{\mathcal{P}}_{2,n}^{\,\mathrm{sq}}(-1,1),\quad A_4=\tfrac{1}{2}\,{\mathcal{P}}_{2,n}^{\,\mathrm{sq}}(-1,-1).
\end{equation}
Substituting~(\ref{eq:identify}) into~(\ref{eq:poly}), this means at
the same time that
\begin{equation}\label{eq:two}
{\mathcal{P}}_n^{\,\mathrm{sq}}(1)=\sum_{p=0}^n(-1)^p\,\sum_{q=0}^p\,\frac{n(n-q-1)!}{{q!(p-q)!}(n-q-p)!}=2.
\end{equation}

\subsection{General case}

We now consider the general case of a square lattice on the torus with $n\times m$ nodes ($n/2\times m/2$ black nodes). In general, the states with a given weight have here a much bigger degeneracy, since the interior of the graph allows for many different configurations which do not affect the boundaries. 
Also here, the highest weight state is unique and in the spin picture
is the one with all spins up.  For large examples, it might seem at
first surprising that the boundary conditions are strong enough to
completely fix the configuration in the interior, but it is easy to
convince oneself by inspecting a small
example
that the interlacing structure of the spins does not allow for any
other configuration in the interior.
This changes once one of the spins on the boundary points down. 
For the next simplest case of $m=4$ for example, each boundary state with weight $n_{max}-1$ has a degeneracy of 3 in the interior. It is obvious that this degeneracy grows quickly with growing $m$ and lower weight. Contrary to the case before, one must always change more than one spin to arrive at a new matching.
On the whole we see that the spin picture is not very well adapted to the general case and the complexity of the derivation used for the strip rises to an unreasonable level for $p > 1$. Also ans\"atze for recursion relations have turned out to be of little use. 
To solve the general problem, other methods might be more appropriate.


\section{Translation to the minimal feedback arc set problem}

The problem of finding a minimal set of arcs in a directed graph upon
the deletion of which the graph becomes acyclic is well studied in
mathematics and computer science~\cite{Karp:1972}. Such a set is
called a minimal \emph{feedback arc set} (\textsc{fas}). There exists
a precise relation between the problem of finding all minimal \textsc{fas} of
a digraph and the problem of identifying all the perfect matchings in
its dual graph. 

Consider a bipartite plane graph $\gra$ with $N$ nodes, and its graph
dual $\gra^\prime$ (where nodes become faces, faces become nodes and edges remain edges) . The dual graph $\gra^\prime$ becomes a digraph if
the edges around a face corresponding to a black node in $\gra$ are
oriented clockwise, while the edges circling a face corresponding to a
white node are oriented counterclockwise. The one--cycles in
$\gra^\prime$ are generated by the plaquettes $\set{p_j}_{j=1}^N$
which correspond to the vertices of $\gra$. Removing the edge $e_{ij}$
shared by the cycles $p_i $ and $p_j$ breaks both cycles. A minimal
\textsc{fas} is obtained taking the collection of $N/2$ edges $\set{e_{i_k
    j_k}}_{k=1}^{N/2}$ shared by disjoint pairs of plaquettes
$p_{i_j}$ and $p_{i_k}$. In the dual graph $\gra$, this corresponds to
selecting a set of edges joining all the nodes and touching all of
them only once. In other words, a minimal \textsc{fas} in $\gra^\prime
$ is a perfect matching in $\gra$.

The situation is different when $\gra $ is embedded on a Riemann
surface of genus $g>0$, because in addition to the one--cycles
generated by the plaquettes, there are $2g$ equivalence classes of
cycles of non--trivial holonomy. In this case, the winding cycles in
$\gra^\prime$ are generated by the zig--zag paths\footnote{A zig--zag
  path is a path which turns alternatingly maximally left and
  maximally right at the vertices.}. It follows that being a perfect
matching in $\gra$ is only a necessary condition for a set of edges to
be a \textsc{fas} in $\gra^\prime$. Let us again restrict ourselves to
the case of $g=1$.  A useful way to represent the partition function
in Eq.~(\ref{eq:torus-partition}) consists in drawing a point in the
$(z,w)$ plane at coordinates $(n_z, n_w)$ for each monomial of the
form $z^{n_z} w^{n_w}$, corresponding to the matchings with weight
$(n_z, n_w)$. In this way, one obtains the region of a $\setZ^2$
lattice delimited by a convex polygon, the \emph{Newton polygon}. We
can hence distinguish between internal and boundary points (or
matchings). It was shown in~\cite{Hanany:2006nm} that a perfect
matching corresponding to an internal point of this Newton polygon is
always a \textsc{fas}. Removing a boundary matching, on the other
hand, always preserves at least one zig--zag path.

\begin{figure}
  \begin{center}
    \subfigure[$(0,1)$ matching]{\includegraphics[width=.4\textwidth]{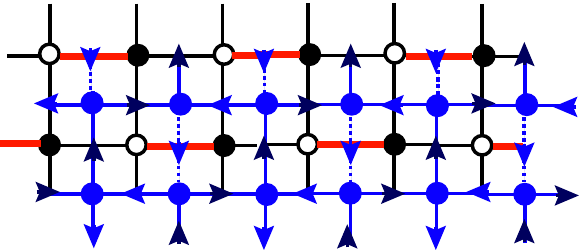}} \hfill
    \subfigure[$(n/2,0)$ matching]{\includegraphics[width=.4\textwidth]{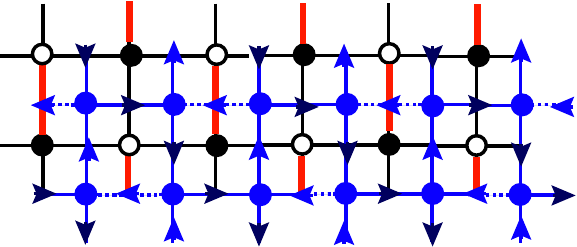}}
    \subfigure[$(0,-1)$ matching]{\includegraphics[width=.4\textwidth]{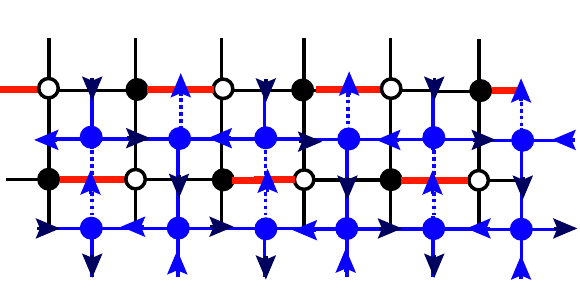}} \hfill
    \subfigure[$(-n/2,0)$ matching]{\includegraphics[width=.4\textwidth]{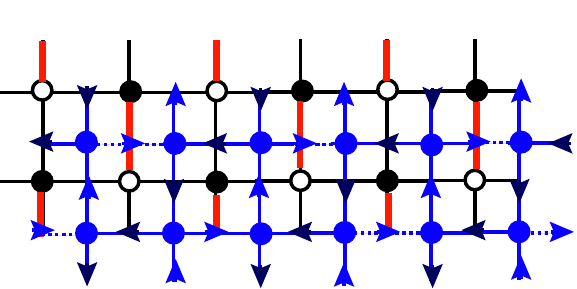}}  
  \end{center}
  \caption{Feedback arc sets on the dual graph to the square strip. The dual graph is represented in light blue. The removal of the boundary matchings (dashed lines) preserves a zig--zag path (dark blue).}
  \label{fig:FAS-square}
\end{figure}

In the case at hand, \emph{i.e.} for the square strip, the polygon is a rhombus
with the four corner vertices at the points $(\pm n/2, 0)$ and
$(0, \pm 1)$. This means that there are four boundary matchings
that will not be feedback arc sets. In fact, removing
the arrows corresponding to these matchings will always preserve
exactly one zig--zag path (see Fig.~\ref{fig:FAS-square}). 
\begin{theo}
  On a $2\times n$ digraph on the torus ($n$ even), there are
  \begin{equation}
    N^{\text{\textsc{fas}}}_n = \left(-1+\sqrt{2}\, \right)^{n} + \left(-1-\sqrt{2}\, \right)^{n} - 2 
  \end{equation}
  minimal feedback arc sets. This number is generated by the function
  \begin{equation}
    \mathcal{F}^{\text{\textsc{fas}}} (s) = \sum_{n=0}^\infty N^{\text{\textsc{fas}}}_n s^n = \frac{2 \left( 1 + \sqrt{2} \right) + \left( \sqrt{2} - 4 \right)s}{\left(s - 1  \right) \left( -1 + \left( -1 + \sqrt{2} \right) s\right)} = \mathcal{F}^{\text{sq}} (s, -1) - \frac{2}{1-s} \, .    
  \end{equation}
\end{theo}
\begin{proof}
  Follows from combining Eq.~(\ref{eq:totalnumber}) and
  Eq.~(\ref{eq:sol}) and subtracting the four boundary
  matchings.
\end{proof}

\section{Discussion and further directions}

In this note, we have derived an explicit formula for the partition function of the dimer model on a $2\times n$ strip of squares or hexagons on the torus. This result was derived via an operator picture, and by mapping the problem to a one--dimensional monomer--dimer system. This result was furthermore translated to the question of finding the minimal feedback arc sets on the dual graph.

An obvious continuation of this work would be the generalization to an $n\times m$ graph. This problem turns out to be much more complicated than the one addressed here and the methods employed here have proven not to be properly adapted to the more general question.


\appendix

\subsection*{Acknowledgements}

We would like to thank Robbert Dijkgraaf for discussions.
Furthermore, we would like to thank the V. Simons Workshop for Mathematics and Physics for hospitality, where this note was concluded.
D.O. is supported in part by \textsc{INFN} and \textsc{MIUR} under
contract 2005-024045 and by the European Community's Human Potential
Program \textsc{MRTN}-\textsc{CT}-2004-005104.  S.R. is supported by
the EC's Marie Curie Research Training Network under the contract
\textsc{MRTN}-\textsc{CT}-2004-512194 "Superstrings".

\newpage
\section{Example: $n=4$ square strip}\label{sec:ex}

\begin{figure}[h!]
  \begin{center}
    \includegraphics[width=.75\linewidth]{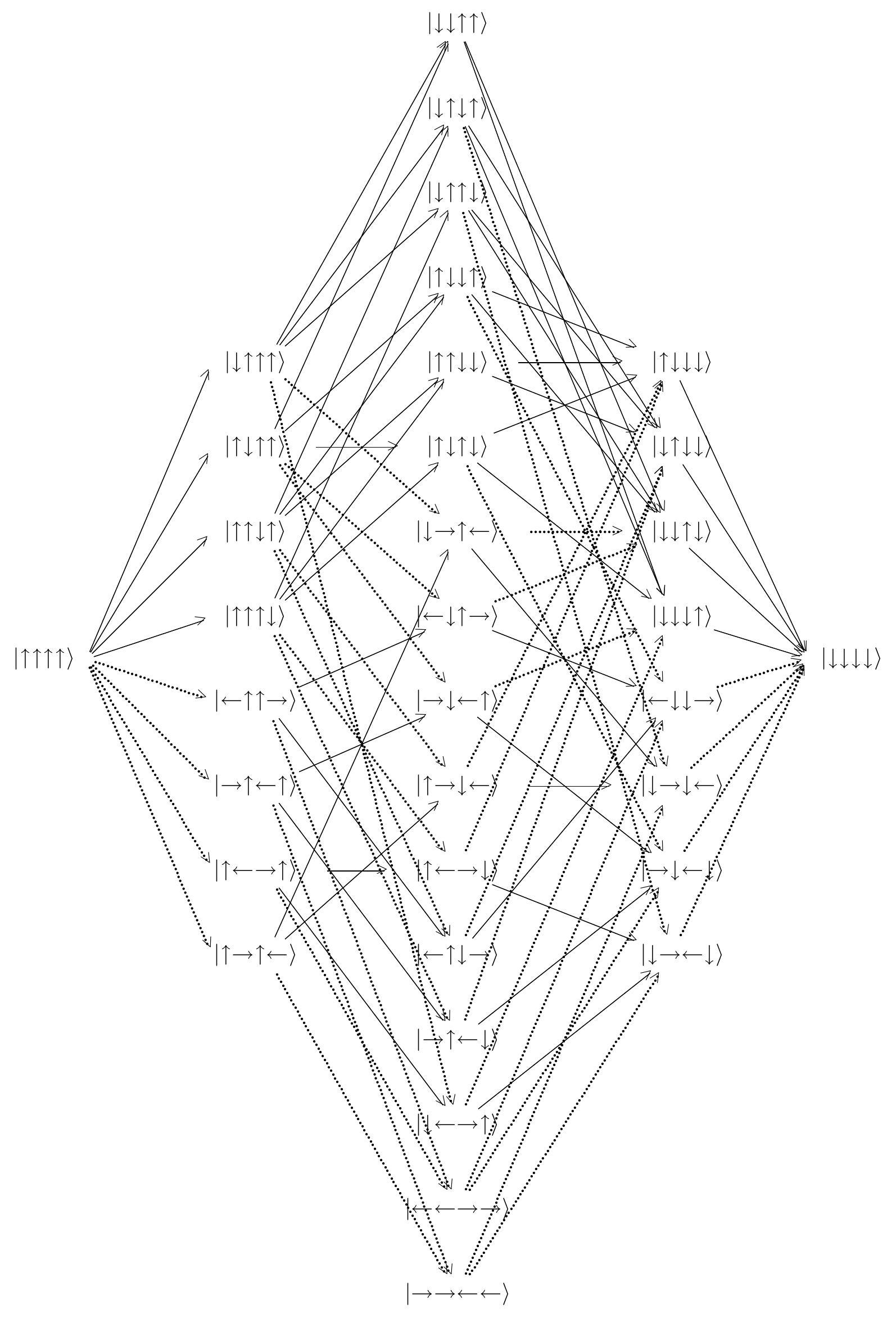}
    \caption{Strip with $n=4$ black nodes}
    \label{fig:sequence}
  \end{center}
\end{figure}

To familiarize ourselves with this operator picture, we consider the
example of $n=4$ black nodes, see Figure \ref{fig:sequence}.  We start
with the highest weight state, which has has weight 2. The solid
arrows in the diagram denote the action by $\di_a$, the dotted arrows
the action by $\di_b$. We find a symmetric sequence structure.


\bibliography{DimerReferences}



\end{document}